\documentclass{article}
\title{Shadow movies not arising from knots}
\author{Daniel Denton\thanks{Research supported by the Neukom Institute for Computation Science, Dartmouth College} \and Peter Doyle}
\date{Version 1.0 dated 17 June 2011}

\usepackage{amsmath}
\usepackage{cancel}
\usepackage{graphicx}

\begin{document}

\maketitle

\begin{abstract}
A shadow diagram is a knot diagram
with under-over information omitted;
a shadow movie is a sequence of shadow diagrams
related by shadow Reidemeister moves.
We show that not every shadow movie arises as the shadow of a
Reidemeister movie,
meaning a sequence of classical knot diagrams related by classical
Reidemeister moves.
This means that in Kaufman's theory of virtual knots,
virtual crossings cannot simply be viewed as classical crossings
where which strand is over has been left `to be determined'.
%
\end{abstract}

\section{Overview}
A shadow diagram is a knot diagram
with under-over information omitted;
a shadow movie is a sequence of shadow diagrams
related by shadow Reidemeister moves.
We will show here that not every shadow movie arises as the shadow of a
Reidemeister movie,
meaning a sequence of classical knot diagrams related by classical
Reidemeister moves.

The question of whether every shadow movie lifts to a Reidemeister movie
arises naturally in the context of virtual knot theory.
(Cf.\ Section \ref{sec:imp} below.)
It also arises very naturally of its own accord:
If a shadow movie can be lifted to a Reidemeister movie, then it can
be lifted to a movie of the unknot, because you can always extend the
shadow movie so that the last frame is an embedded circle.
Now, it is a staple of introductory courses on knot theory to show that any shadow diagram arises as a shadow of the unknot:
Just imagine laying a string down along the curve, starting at any point you
please,
until you get back to where you started.
If a shadow movie can be morphed so that a simple point $p$ of the curve remains
simple throughout the movie, then resolving each frame to the unknot using
$p$ as the starting point yields a lift to a Reidemeister movie of the unknot.
Put another way, in the category of `long knots',
all shadow movies lift to movies of the `long unknot'.
But what if all points of the knot are necessary `implicated'?
Can the shadow movie still be lifted?

\begin{figure}[htb]
        \center{\includegraphics[scale=0.18]
        {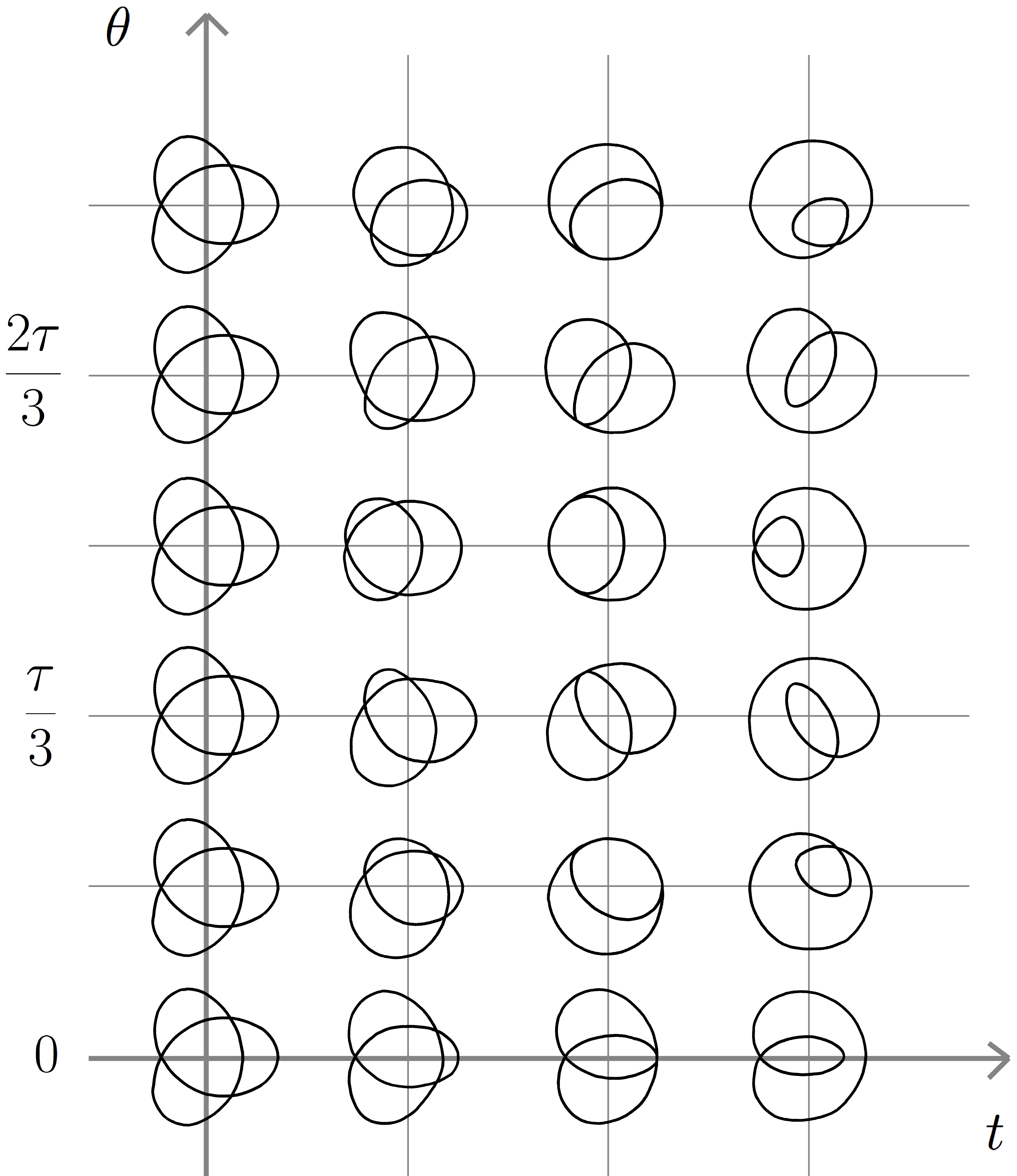}}
        \caption{An unliftable 2-parameter family of shadows.
We parametrize a family of perturbations
of the standard trefoil shadow diagram
by two parameter $(t, \theta)$,
where $t$ is a time parameter and $\theta$ a direction parameter.
Fixing $\theta$ to be $0$, $\tau/3$ or $2\tau/3$ and varying $t$, we
see each of the 3 possible type II shadow moves.
No designation of the crossings is simultaneously compatible with all
3 moves, so this 2-parameter family doesn't lift.}
\label{two_parameter}
\end{figure}

Or consider this:
One readily shows that a two-parameter family of shadows
cannot always be lifted
(see Figure~\ref{two_parameter}).
But playing around with small examples seems to allow the possibility
that one-parameter families might always be liftable.

Our example of an unliftable shadow movie was found by means of a
computer search.
The computer generated shadow movies by
walking at random through the space of shadow diagrams
(or rather, the space of associated chord diagrams) while keeping track of which assignments of signs to the crossings
were compatible with the movie so far.
Eventually it found cases of movies where the set of compatible signs became empty,
meaning that the shadow movie was unliftable.
The example presented here was found by massaging one of
the unliftable movies that the computer found.

As far as we know, the result reported here is new.
On the other hand, the question we answer here seems so natural
that we presume others must have considered it.
And the example is simple enough that it---or perhaps another,
even simpler example---could have been found by hand.

\section{Preliminaries}

Let a \emph{shadow diagram} be a classical knot diagram where, at each crossing, the information
about which arc passes on top and which one passes underneath has been removed.
Thus a shadow diagram is an immersion of the circle into the plane in such a 
way that all the intersections are transverse crossings of exactly two arcs.

There is a natural mapping of every classical knot diagram to a shadow diagram
obtained by removing the under-over information at each crossing.
Since there are $2$ ways in which each 
crossing can be arranged, a shadow diagram with $n$ crossings arises from $2^n$ 
different classical diagrams.

\begin{figure}[ht]
\begin{minipage}[b]{0.45\linewidth}
\centering
\includegraphics[scale=0.5]{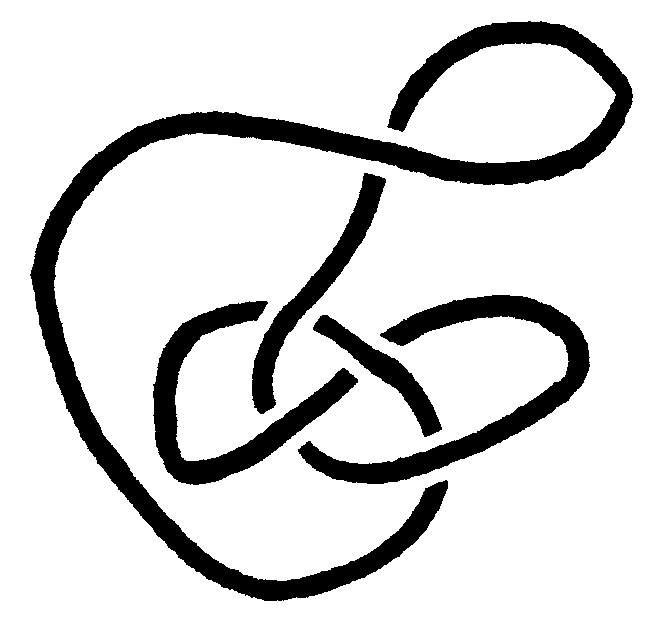}
\caption{A Classical Knot Diagram}
\end{minipage}
\hspace{0.1\linewidth}
\begin{minipage}[b]{0.45\linewidth}
\centering
\includegraphics[scale=0.5]{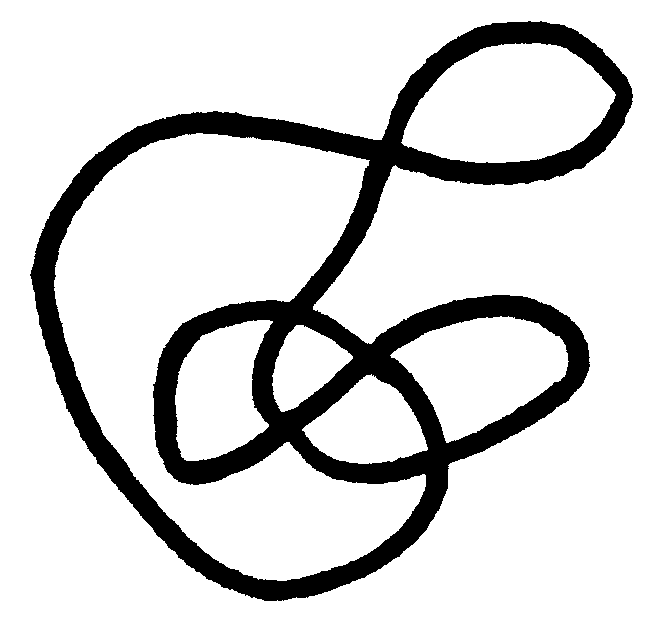}
\caption{A Shadow Knot Diagram}
\end{minipage}
\end{figure}

The sign of a crossing is a way of indicating which of the two arrangements
or resolutions a crossing takes on (in a classical knot diagram).  If we assign an orientation
(or direction of travel) to the diagram, then the sign is positive if the angle between the direction 
of travel on the over-strand and the direction of travel on the under-strand is between $0$ and $\tau/2$.
Otherwise the sign is negative.

\begin{figure}[ht]
\begin{minipage}[b]{0.45\linewidth}
\centering
\includegraphics[scale=0.5]{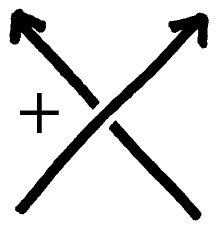}
\caption{A Crossing With Positive Sign}
\end{minipage}
\hspace{0.1\linewidth}
\begin{minipage}[b]{0.45\linewidth}
\centering
\includegraphics[scale=0.5]{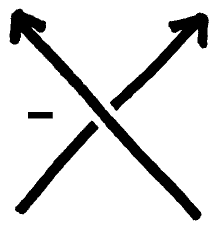}
\caption{A Crossing With Negative Sign}
\end{minipage}
\end{figure}

Conveniently, in the case of single component knots the signs of the crossings are not affected by which orientation we pick for the curve
(since reversing the orientation reverses the direction of travel along both strands simultaneously).
Thus we can talk about the signs of the crossings of a classical knot diagram even though
we have not indicated an orientation.  

Just as Reidemeister moves can be used to operate on classical knot diagrams, we can 
define corresponding `shadow Reidemeister moves' to operate on shadow diagrams.
An important difference is that, while we use equivalence classes of classical diagrams 
under Reidemeister moves to define the type of a classical knot, there is only the one
trivial equivalence class of shadow diagrams under shadow Reidemeister moves,
as all shadow diagrams are equivalent to a diagram with no crossings through a sequence of shadow moves.
(Cf.\ Fenn and Kauffman \cite{fenn_vktunsolvedproblems}.)

\begin{figure}[ht]
\begin{minipage}[b]{0.45\linewidth}
\centering
\includegraphics[scale=0.333333]{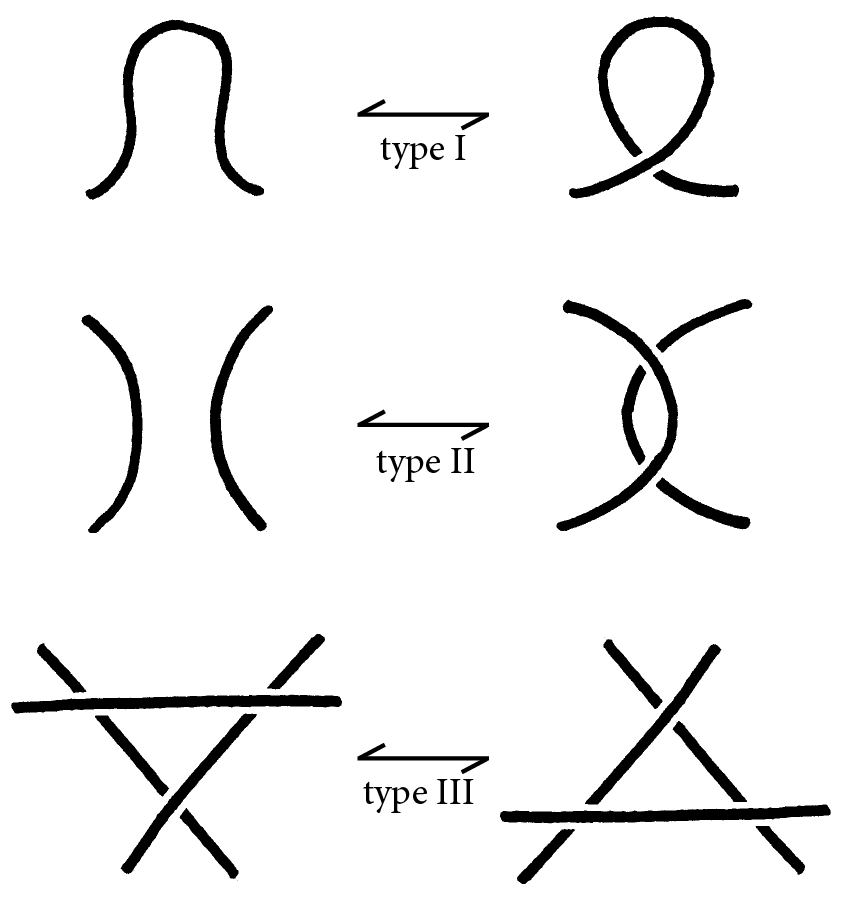}
\caption{Classical Reidemeister Moves}
\end{minipage}
\hspace{0.1\linewidth}
\begin{minipage}[b]{0.45\linewidth}
\centering
\includegraphics[scale=0.333333]{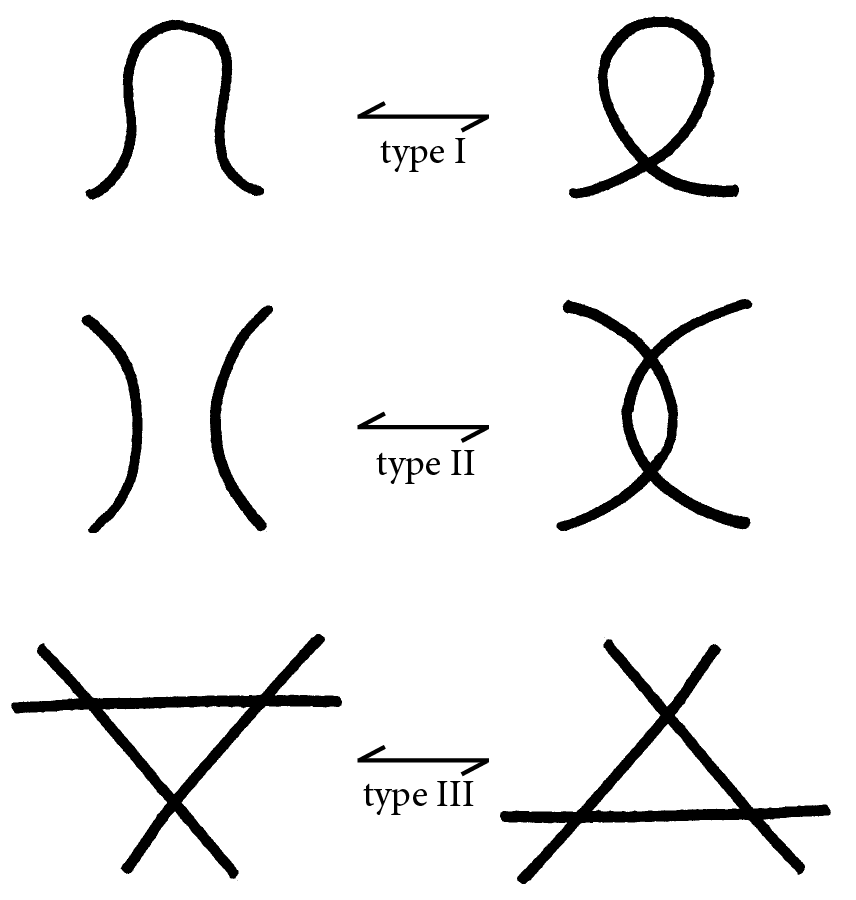}
\caption{Shadow Reidemeister Moves}
\end{minipage}
\end{figure}

We use the term \emph{Reidemeister movie}
to refer to a sequence of classical knot diagrams
connected by Reidemeister moves.  (Thus two classical knot diagrams represent the 
same knot if and only if they are connected by a Reidemeister movie.)
Analogously, by a \emph{shadow movie}
we mean a sequence of shadow diagrams
connected by shadow Reidemeister moves.

We will show that not every valid shadow movie can be lifted to a valid Reidemeister movie
through some assignment of crossing signs by presenting an example of an unliftable shadow movie.  

\section{The example}
While we could present this example as a sequence of individual Reidemeister moves,
doing so would be very long and cumbersome,
so we permit ourselves the following 
two simplifications.

\subsection*{Type $0$ Moves on the Sphere}

When we have a loop on the outside edge of our diagram, we can redraw this loop 
around the other side of the diagram since we can clearly pull the entire loop over (or under)
the rest of the diagram without affecting the sign constraints on any of the already existing 
crossings.  We call this a type $0$ move on the sphere (abbreviated $0S^2$) since, were we to draw our diagram
on a sphere instead of the plane, we could accomplish this without any Reidemeister moves by
simply swinging the loop around across the far side of the sphere.

\begin{figure}[htb]
        \center{\includegraphics[width=\textwidth]
        {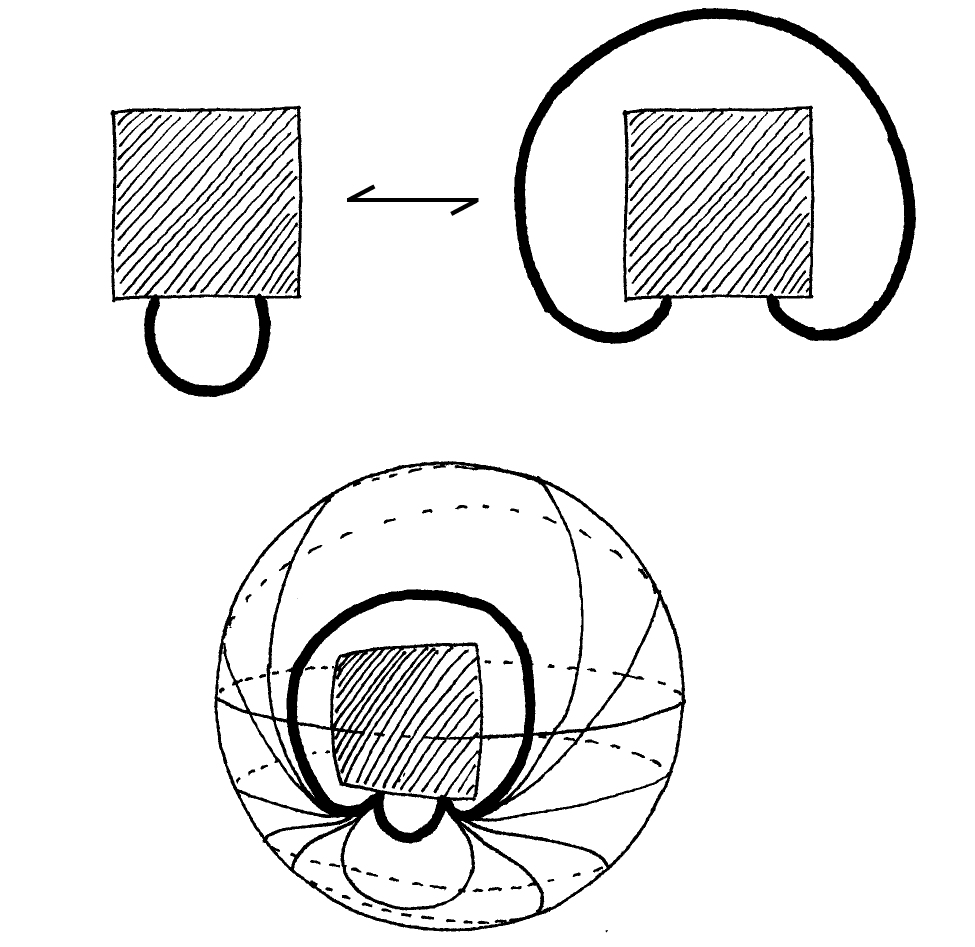}}
        \caption{A Type 0 Move on the Sphere}
\end{figure}

\subsection*{Ear-rolling Move}

Whenever we have a ear on our diagram (a loop as created by a type I Reidemeister
move) along some section of the curve, we can roll the ear along the the curve to another location
where it doesn't cross with any other strands.  Once again, this procedure does not
affect the constraints placed on any of the of the already existing crossings.

\begin{figure}[htb]
        \center{\includegraphics[width=\textwidth]
        {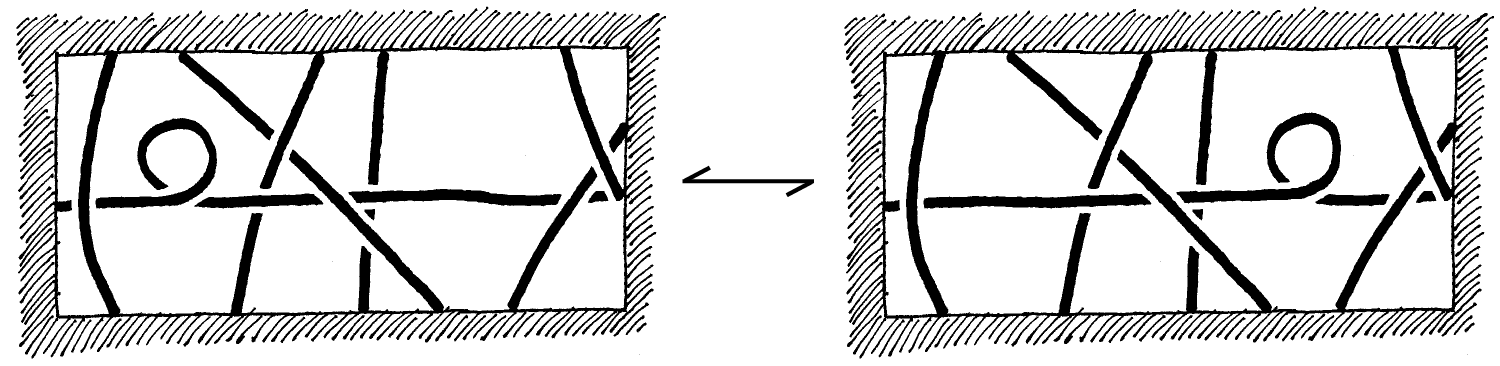}}
        \caption{An `Ear-rolling' Move}
\end{figure}

\subsection*{Example}

We start with the unknot and perform a type II shadow Reidemeister move.  Because type 
two moves can only introduce or remove crossings with opposite signs, this introduces 
a constraint on the crossings' signs.  We denote the signs of the two crossings as $a$ and
$-a$ to show that, whichever sign the first is, it is the opposite of the second.

\begin{figure}[htb]
        \center{\includegraphics[scale=0.5]
        {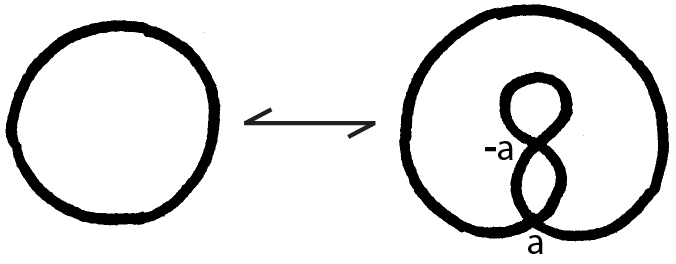}}
\end{figure}

Through a series of ear-rolling and type $0S^2$ moves we then rearrange the shadow diagram
without affecting the crossing constraints.

\begin{figure}[htb]
        \center{\includegraphics[scale=0.5]
        {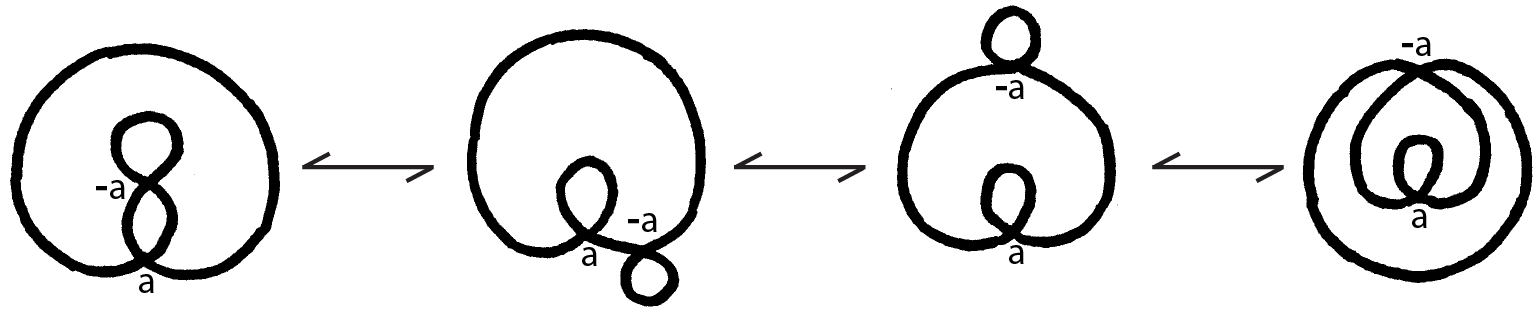}}
\end{figure}

\newpage
A second type II move introduces two more crossings with signs $b$ and $-b$.

\begin{figure}[htb]
        \center{\includegraphics[scale=0.5]
        {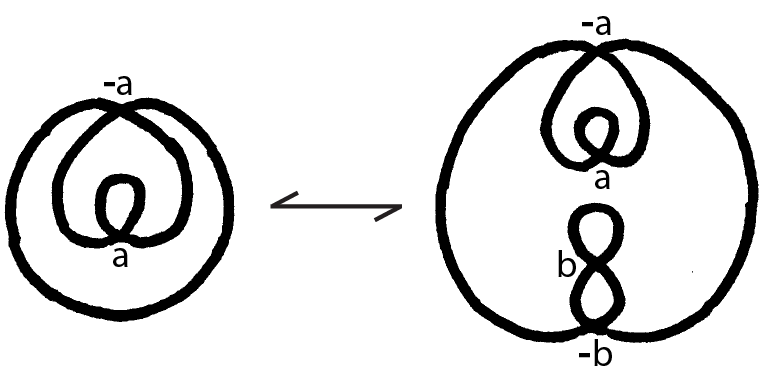}}
\end{figure}

This is followed by more ear-rolling and type $0S^2$ moves.

\begin{figure}[htb]
        \center{\includegraphics[scale=0.5]
        {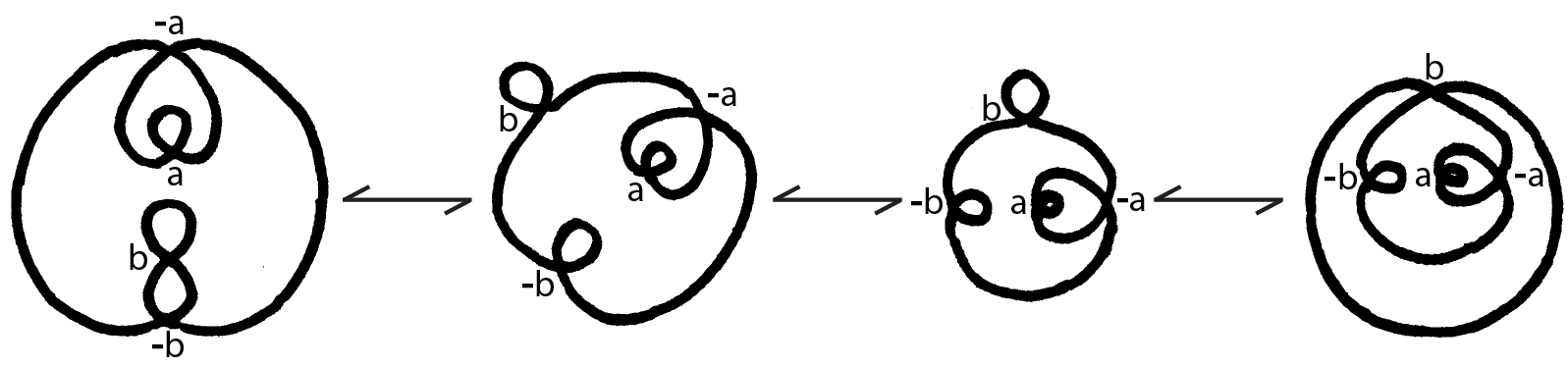}}
\end{figure}

Finally, one more type II move introduces a fifth and sixth crossing with signs $c$ and $-c$ respectively.

\begin{figure}[htb]
        \center{\includegraphics[scale=0.5]
        {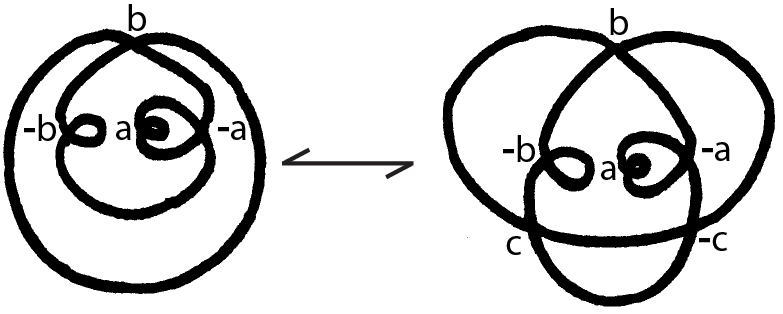}}
\end{figure}

We then perform one more type $0S^2$ move and one more ear-rolling move.

\begin{figure}[htb]
        \center{\includegraphics[scale=0.5]
        {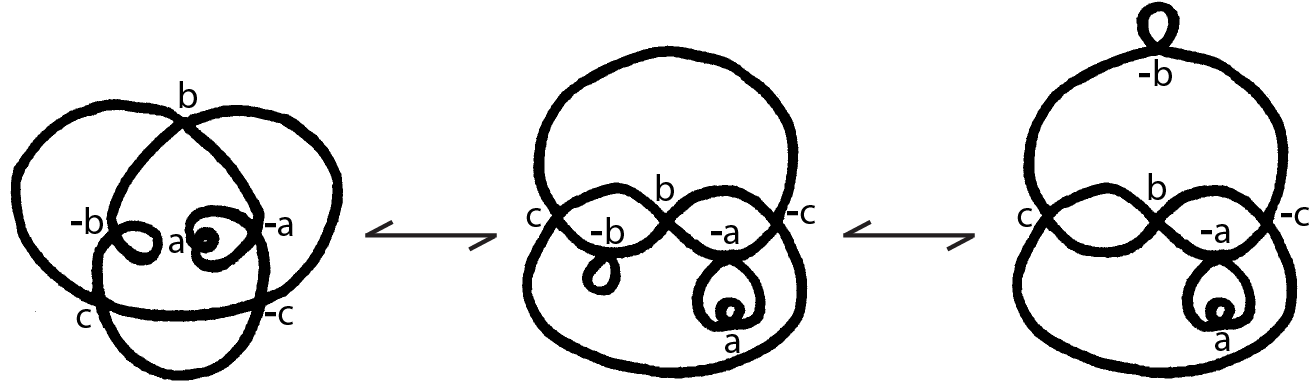}}
\end{figure}

\newpage
The next move is a type III Reidemeister move, collapsing a triangle whose corner crossings
have the signs $-a$, $b$, and $-c$.  Because of the directions of the orientations of the segments
between the $3$ crossings, one of the following three relations must hold.
This particular set of
constraints will not turn out to be needed in our proof that
this shadow movie unliftable (which suggests 
that this may not be a minimal such example), but we mention it anyway.

\begin{figure}[ht]
\begin{minipage}[b]{0.45\linewidth}
\centering
\includegraphics[scale=0.5] {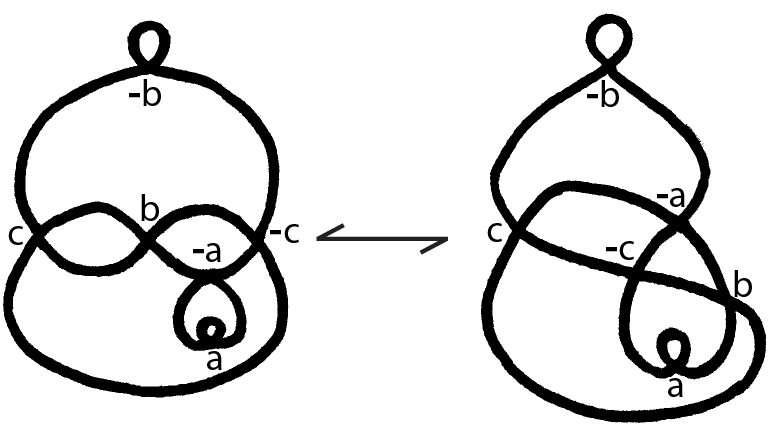}
\end{minipage}
\hspace{0.1\linewidth}
\begin{minipage}[b]{0.45\linewidth}
\centering
\begin{eqnarray*}
 & -a=b=-c\\
 & \text{or}\\
 & -a=b=c\\
 & \text{or}\\
 & -a=-b=-c
\end{eqnarray*}
\end{minipage}
\end{figure}

This is followed by another type III Reidemeister move imposing another three possible relations.

\begin{figure}[ht]
\begin{minipage}[b]{0.45\linewidth}
\centering
\includegraphics[scale=0.5] {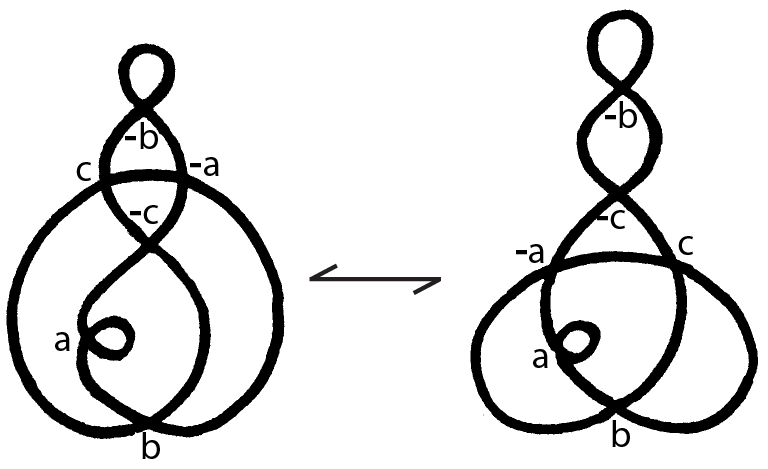}
\end{minipage}
\hspace{0.1\linewidth}
\begin{minipage}[b]{0.45\linewidth}
\centering
\begin{eqnarray*}
 & \cancel{-c=c=-a}\\
 & \text{or}\\
 & \cancel{-c=c=a}\\
 & \text{or}\\
 & -c=-c=-a
\end{eqnarray*}
\end{minipage}
\end{figure}

Because the first two of these relations are self-contradicting, we are constrained to have $a=c$.

\newpage
After a third type III Reidemeister move we get three more possible relations.

\begin{figure}[ht]
\begin{minipage}[b]{0.45\linewidth}
\centering
\includegraphics[scale=0.5] {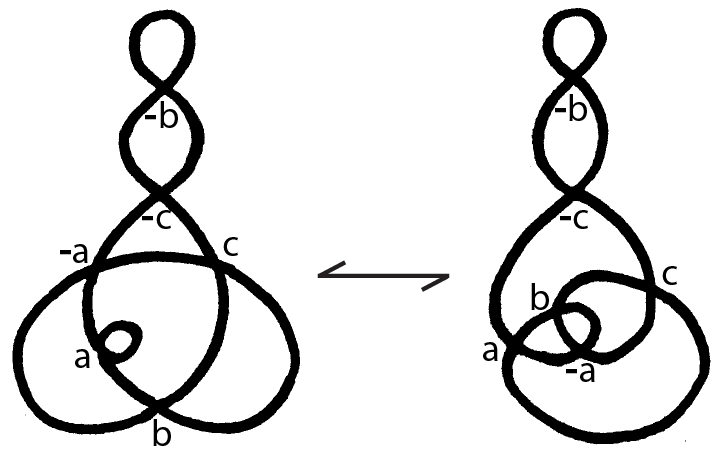}
\end{minipage}
\hspace{0.1\linewidth}
\begin{minipage}[b]{0.45\linewidth}
\centering
\begin{eqnarray*}
 & \cancel{a=-a=b}\\
 & \text{or}\\
 & \cancel{a=-a=-b}\\
 & \text{or}\\
 & a=a=b
\end{eqnarray*}
\end{minipage}
\end{figure}

Once again, two of these relations are self-contradicting leaving us with $a=b$.  Thus we now 
have the signs of all $6$ crossings related through the relation $a=b=c$.

The example finishes with a type II shadow Reidemeister move removing two crossings.
These crossings must have opposite signs in order to carry out the move, so $b=-c$.

\begin{figure}[htb]
        \center{\includegraphics[scale=0.5]
        {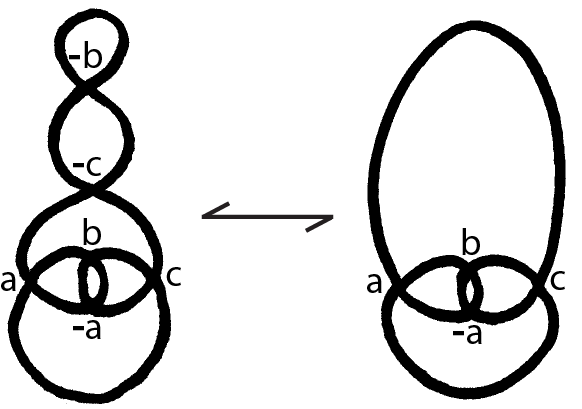}}
\end{figure}

This, however, contradicts the already imposed constraints.  Therefore, there is no 
resolution of crossing signs which is consistent with this sequence of Reidemeister moves,
and hence this is an example of a shadow movie which cannot be lifted to a valid Reidemeister movie.

\section{Implications} \label{sec:imp}

One of the fundamental results of virtual knot theory is Kauffman's theorem
that any two classical knot diagrams which are equivalent under extended
virtual Reidemeister moves are also equivalent under classical Reidemeister moves.
(Cf.\ \cite{kauffman_vkt}.)

This theorem would follow immediately if it were true that any virtual knot movie could be lifted to a Reidemeister movie through some choice of over and under strands at the virtual crossing.
In fact, we would have the stronger result that the minimum number of moves to pass from one classical knot diagram to another equivalent diagram does not decrease when you allow extended virtual Reidemeister moves along with the classical
moves.
Our example lays to rest the naive hope that virtual crossings are just classical crossings whose signs are `to be determined',
because we can't even lift virtual movies of the unknot.
However, we still don't know any example where allowing virtual moves decreases
the number of moves needed to get from one classical knot projection to another.

\bibliographystyle{plain}
\bibliography{shadow}

\begin{thebibliography}{1}

\bibitem{fenn_vktunsolvedproblems}
Roger Fenn, Louis~H. Kauffman, and Vassily~O. Manturov.
\newblock Virtual knot theory - unsolved problems.
\newblock {\em Fundamenta Mathematicae}, 188:293--323, 2005.

\bibitem{kauffman_vkt}
Louis~H. Kauffman.
\newblock Virtual knot theory.
\newblock {\em European Journal of Combinatorics}, 20:663--691, 1999.

\end{thebibliography}

\end{document}